\documentclass[12pt]{amsart}
\usepackage{graphicx}
\usepackage{amssymb}
\usepackage{amscd}
\usepackage{hyperref}
\usepackage{verbatim}
\usepackage{amsmath}
\usepackage{marvosym}
\newtheorem{thm}{Theorem}
\newtheorem{theorem}{Theorem}[section]
\newtheorem{corollary}[theorem]{Corollary}
\newtheorem{lemma}[theorem]{Lemma}
\newtheorem{prop}[theorem]{Proposition}

\theoremstyle{remark}
\newtheorem{rem}[theorem]{\bf Remark}

\theoremstyle{definition}
\newtheorem{definition}[theorem]{Definition}

\theoremstyle{definition}
\newtheorem*{ack}{Acknowledgments}

\theoremstyle{remark}

\newcommand{\dbar}{\bar\partial}

\newcommand{\dint}{d_\square} %or {d_{\dbar}}
\newcommand{\dom}{{\rm dom\,}}
\newcommand{\domo}{{\rm dom}_0\,}
\newcommand{\domc}{{\rm dom}_c\,}

\newcommand{\supp}{\, {\rm supp}\,}
\newcommand{\bid}{\, {\bf 1}}

\newcommand{\mbar}{M\llap{\raise 1.2ex\hbox{---}}}
\newcommand{\dirac}{{\partial}\llap{\raise 0.1ex\hbox{/}}}

\newcommand{\lpq}{\Lambda^{p,q}}

\newcommand{\cbpq}{C^\infty(\mbar,\lpq)}
\newcommand{\ccpq}{C^\infty_c(\mbar,\lpq)}
\newcommand{\lzpq}{L^2(M,\Lambda^{p,q})}
\newcommand{\lupq}{L^\infty(M,\Lambda^{p,q})}

\newcommand{\esssup}{\, {\rm esssup}}
\newcommand{\vol}{\, {\rm vol}\,}
\newcommand{\ZZ}{\mathbb{Z}}
\newcommand{\RR}{\mathbb{R}}
\newcommand{\NN}{\mathbb{N}}
\newcommand{\cH}{\mathcal{H}}
\newcommand{\CC}{\mathbb{C}}

\newcommand{\Hmm}[1]{\leavevmode{\marginpar{\tiny%
$\hbox to 0mm{\hspace*{-0.5mm}$\leftarrow$\hss}%
\vcenter{\vrule depth 0.1mm height 0.1mm width \the\marginparwidth}%
\hbox to 0mm{\hss$\rightarrow$\hspace*{-0.5mm}}$\\\relax\raggedright #1}}}
\begin{document}
\title[Generalized eigenforms]{Essential self-adjointness, generalized eigenforms, and spectra for the $\dbar$-Neumann problem on $G$-manifolds.}
\author{Joe J Perez\ \ }
\address{Fakult\"at f\"ur Mathematik\\
Universit\"at Wien\\
Vienna, Austria}
\email{joe\_j\_perez@yahoo.com}
\thanks{JJP is supported by FWF grants P19667 and I382}

\author{\ \ Peter  Stollmann}
\address{Faculty of Mathematics\\
Chemnitz University of Technology\\
Chemnitz, Germany}
\email{peter.stollmann@mathematik.tu-chemnitz.de}
\thanks{PS is partly supported by DFG}

\begin{abstract} Let $M$ be a strongly pseudoconvex complex manifold which is
also the total space of a principal $G$--bundle with $G$ a Lie group and compact orbit space $\mbar/G$. Here we investigate the $\dbar$-Neumann
Laplacian $\square$ on $M$. We show that it is essentially self-adjoint on its
restriction to compactly supported smooth forms. Moreover we relate its spectrum to the existence of generalized eigenforms: an energy belongs to $\sigma(\square)$ if there is a subexponentially bounded generalized
eigenform for this energy. Vice versa, there is an expansion in terms of
these well--behaved eigenforms so that, spectrally, almost every energy comes with such a generalized
eigenform.
\end{abstract}

%\date{\today}
%\addcontentsline{toc}{chapter}{}

\maketitle
\tableofcontents

%%%%%%%%%%%%%%%%%%%%
\section{Introduction}
%%%%%%%%%

The approach to the theory of several complex variables via partial differential equations involves the analysis of a self-adjoint boundary value problem for an operator $\square$ similar to the Hodge Laplacian. This problem, called the $\dbar$-Neumann problem, is  the subject of this article and we will give a brief description here.

We will assume that $M$ is a complex manifold, $n=\dim_{\mathbb C} M$, with
smooth boundary $bM$ such that $\mbar = M\cup bM$. Assume further that $\mbar$
is strictly contained in a slightly larger complex manifold $\widetilde{M}$ of the same
dimension. For any integers $p,q$ with $1\leq p,q\leq n$ denote by
$C^\infty(M,\Lambda^{p,q})$ the space of all $C^\infty$ forms of type $(p,q)$ on
$M$, {\it i.e.}\ the forms which can be written in local complex coordinates
$(z^1, z^2,\dots,z^n)$ as
\begin{equation}\label{pqform}\phi=\sum_{|I|=p,|J|=q}\phi_{I,J}\ dz^I\wedge
d\bar z^J \end{equation}
\noindent
where $dz^I=dz^{i_1}\wedge\dots\wedge dz^{i_p}$, $d\bar z^J=d\bar
z^{j_1}\wedge\dots\wedge d\bar z^{j_q}$, $I=(i_1,\dots,i_p)$, $J=(j_1,\dots,
j_q)$, $i_1<\dots<i_p$, $j_1<\dots<j_q$, with the $\phi_{I,J}$ smooth
functions in local coordinates.  For such a form $\phi$, the value of the
antiholomorphic exterior derivative $\dbar\phi$ is
\[\dbar\phi=\sum_{|I|=p,|J|=q} \sum_{k=1}^n
\frac{\partial\phi_{I,J}}{\partial\bar z^k}\ d\bar z^k\wedge dz^I\wedge d\bar
z^J\]
\noindent
so $\dbar = \dbar|_{p,q}$ defines a linear map
$\dbar:C^\infty(M,\Lambda^{p,q})\to C^\infty(M,\Lambda^{p,q+1})$.
With respect to a smooth measure $\mu$ on $M$ and a smoothly varying Hermitian
structure in the fibers of the tangent bundle, define the spaces
$L^2(M,\Lambda^{p,q})$. Let us consider $\dbar$ as the maximal operator in $L^2$
and let $\dbar^\ast$ be its Hilbert space adjoint operator (this involves the
introduction of boundary conditions). Define the nonnegative form
\begin{equation}\label{qform}Q(\phi,\psi)
=\langle\dbar\phi,\dbar\psi\rangle_{L^2(M,\Lambda^{p,q+1})} +
\langle\dbar^*\phi,\dbar^*\psi\rangle_{L^2(M,\Lambda^{p,q-1})},\end{equation}
\noindent
with domain $\dom(Q)=\dom(Q^{p,q})\subset L^2(M,\Lambda^{p,q})$ and denote the
associated self-adjoint operator in $L^2(M,\Lambda^{p,q})$ by
$$\square=\square_{p,q}=\dbar^\ast\dbar + \dbar\dbar^\ast,
$$
using $+$ for the \emph{form sum} of two
self-adjoint operators; see \cite{F}.
The Laplacian is elliptic but its natural boundary conditions are not coercive,
thus, in the interior of $M$, the operator gains two degrees in the
Sobolev scale, as a second-order operator, while in neighborhoods of the boundary, it
gains less. The gain at the boundary depends on the geometry of the boundary,
and the best such situation is that in which the boundary is {\it strongly
pseudoconvex}. In that case, the operator gains one
degree on the Sobolev scale in neighborhoods of $bM$ and so global estimates including both interior and
boundary neighborhoods gain only one degree.

More generally, one says that the Laplacian satisfies a {\it pseudolocal
estimate with gain $\epsilon>0$ in} $L^2(M,\Lambda^{p,q})$ in the following situation.

 If $U\subset\mbar$ is a neighborhood with compact closure, $\zeta, \zeta'\in
C^{\infty}_{c}(U)$ for which $\zeta'|_{{\rm supp}(\zeta)}=1$, and $\alpha|_U\in
H^{s}(U,\Lambda^{p,q})$, then $\zeta(\square +\bid )^{-1}\alpha\in
H^{s+\epsilon}(\mbar,\Lambda^{p,q})$ and there exists a constant $C_{\zeta, \zeta'}>0$ such that
\begin{equation}\label{prima1}\|\zeta (\square
+\bid)^{-1}\alpha\|_{H^{s+\epsilon}(M,\Lambda^{p,q})}\le C_{\zeta, \zeta'}
(\|\zeta'\alpha\|_{H^s(M,\Lambda^{p,q})}+\|\alpha\|_{L^2(M,\Lambda^{p,q})})
\end{equation}
\noindent
uniformly for all $\alpha$ satisfying the assumption. See \cite{K1, K2, KN, FK, E} for these results.

Mostly for the simplicity that a group symmetry implies, let us in this paper that the manifolds in consideration satisfy the following requirements.

\begin{definition} \label{assumption} We will say that $M$ satisfies
\emph{Assumption} {\bf (A)} if the following hold.
First, assume that $M$ is a complex manifold which is also the total space of a principal $G$-bundle
 with $G$ a Lie group acting by holomorphic transformations and with compact orbit space $\mbar/G$:
\[G \longrightarrow M \longrightarrow X.\]
\noindent
Also assume also that $M$ has a strongly pseudoconvex boundary as above so that $\square=\square_{p,q}$ satisfies a pseudolocal estimate with gain $\epsilon=1$ in $L^2(M,\Lambda^{p,q})$ with $q>0$. \end{definition}

Though our results hold in substantially greater generality, which we will indicate where we feel necessary, we keep our setting as above, with exact invariances. We note that in the case in which $G$ is unimodular, there is a good generalized Fredholm theory for the $\square$ as well as generalized Paley-Wiener theorems for $G$-bundles which together provide an effective framework for understanding the solvability of equations involving $\square$. These are worked out and applied in \cite{P1, P2, DSP}. In \cite{P3}, the unimodularity condition is dropped, as in our setting here.

Bundles constructed in \cite{GHS, HHK} are examples of manifolds satisfying our assumptions with $\epsilon=1$ for all $q>0$.

We will in this article be concerned with the following fundamental properties of the operator $\square$.

\begin{thm} Assume {\rm \textbf{(A)}} from \ref{assumption}. Then $\square$ is essentially self-adjoint on $C^\infty_c(\mbar,\Lambda^{p,q})\cap \dom\square$.\end{thm}

This type of result is very common for many natural partial differential operators on manifolds without boundary. It is important because it provides that there is only one way to extend the operator from a domain consisting of smooth, compactly supported forms to a self-adjoint operator. The case at hand is more complicated due to the boundary, which moreover plays an important role and comes with noncoercive boundary conditions. We prove Theorem 1 in Sect.\ \ref{essential} by a cutoff procedure that requires taking the boundary condition into account. We borrow from \cite{BS} and from discussions of the first-named author with E.\ Straube.

The reader will notice that we state this theorem first and prove it last. The reason for this is that we base the following two results on quadratic form methods for which we do not need the more precise description of the domain of the operator.

\begin{thm}\label{schnol-thm} \emph{\bf (Schnol-type theorem)}  Assume {\rm
\textbf{(A)}} from \ref{assumption}. The existence
of a generalized eigenform for $\square$ with
eigenvalue $\lambda$ satisfying certain growth
conditions implies that $\lambda\in\sigma(\square)$.
\end{thm}

This type of result is often called Schnol's theorem in the literature. The
precise statements are Theorem \ref{half} and Corollary \ref{sub} below.
Actually, the original result of Schnol's paper \cite{Sh} is an equivalence, so
Theorems \ref{schnol-thm} and the following Theorem \ref{thm1.3} together give results reminiscent of Schnol's theorem from
Schr\"odinger operator theory; see
\cite{Sh} and the discussion in \cite{BdMLS, LSV1} for a list of references
and recent results in the Dirichlet form context.
See also \cite{Shu1,Shu2} for results on general elliptic operators on sections
of vector bundles over complete manifolds.

\begin{thm}\label{thm1.3} \emph{\bf (Eigenfunction expansion)} Assume {\rm
\textbf{(A)}} from \ref{assumption}. Let
$\omega\in L^2(M,\RR)$ with $\omega^{-1}\ge 1$. Then, for spectrally a.e.\
$\lambda\in\sigma(\square)$ there is a generalized eigenform
$\varepsilon_\lambda$ for $\square$ with eigenvalue $\lambda$ so that
$\omega\varepsilon_\lambda\in \lzpq$. \end{thm}

For the proof of Theorem \ref{schnol-thm} we follow the strategy from \cite{BdMLS}, see also \cite{LSV1}: starting from the well behaved
generalized eigenform $u$ we construct a singular sequence $u_k=\eta_k u$ for the form $Q$ of $\square$. The cutoff functions have to be such that the product $\eta_k u$ belongs to the domain of the form $Q$. That is achieved by using the intrinsic metric of $\square$ to define $\eta_k$. The intrinsic metric for $\square$ was introduced in our previous work \cite{PS} and turned out to be useful in estimating the heat kernel of $\square$. Here we provide some more results and a useful characterization of the intrinsic metric in Section \ref{metric_sec}. That the cutoff does provide a singular sequence is a consequence of a Caccioppoli type inequality, which is the subject of Section \ref{Sec_4}. In Section \ref{Sec_5}
 we prove two variants of Theorem \ref{schnol-thm}, making precise what  ``certain growth
conditions'' means.

Expansion in generalized eigenelements is typically based on strong compactness properties. Here we use the method developed in \cite{BdMS} for Dirichlet forms, based on an abstract result from \cite{PSW}. The main input is from \cite{PS}, where we showed ultracontractivity of the heat semigroup corresponding to $\square$; we also refer to this paper for more pointers to related literature.

%%%%%%%%%%%%%%%%%%%%%%%%%%%
\section{Preliminaries}
%%%%%%%%%%%%%%%%%%%%%%%%%%%
\label{Prel}
We will have to
describe smoothness of functions, forms, and sections of vector bundles using
$G$-invariant Sobolev spaces which
we define here.

We denote by $C^\infty(M,\lpq)$ the
space of smooth $(p,q)$-forms on $M$, by $\cbpq$ the subspace of those forms
that can be smoothly extended to $\mbar$ and by $\ccpq$ the subspace of the
latter, consisting of those smooth forms with compact support. Given any
$G$-invariant, pointwise Hermitian structure
\[C^\infty(\mbar,\Lambda^{p,q})\ni u,v \longmapsto \langle u(x) ,v(x)
\rangle_{\Lambda_x^{p,q}} \in\CC, \quad (x\in \mbar),\]
\noindent
and its volume form $\mu$, we define the $L^p$-spaces $L^p(M,\Lambda^{q,r})$ of forms, for
$1\le p\le\infty$, as those forms $u$, for which the norm
\[\|u\|_{L^p(M,\Lambda^{q,r})} = \left[\int_M \langle
u,u\rangle_{\Lambda^{q,r}}^{p/2}\, d\mu\right]^{1/p} .\]
\noindent
is finite, with the obvious modification for $p=\infty$.
We will sometimes abbreviate $\langle u,u\rangle_{\Lambda^{q,r}}$ by writing
$|u|_{\Lambda^{q,r}}^2$ instead. Also, we will write
$\langle\cdot,\cdot\rangle_{\Lambda^p}$ to mean the Hermitian structure on
$\mathbb C\otimes\Lambda^p =\oplus_S\Lambda^{q,r}$ with $S=\{(q,r)\mid q+r=p\}$
as well as the Riemannian metric on $\Lambda^1$ associated, see
\cite[\S3.2]{PS}.

As we have a manifold with bounded geometry, there exist partitions of unity
with bounded multiplicity and derivatives, \cite{Gro, Gro1, Ko1, Ko2, Schick-96, Shu2}
and, by differentiating componentwise with respect to local geodesic coordinates,
we may assemble $G$-invariant integer Sobolev spaces $H^s(M,\Lambda^{p,q})$, for
$s=0,1,2,\dots$.

Because $\mbar/G$ is compact, the spaces
$H^s(M,\Lambda^{p,q})$ do not depend on the choice of an invariant Hermitian structure on
$\Lambda^{p,q}$.
The usual duality relations for $L^p$ spaces hold (polarizing
the above norm) as well as the Sobolev lemma, {\it
etc.} Background on this is provided in \cite{GKS}.
We will also need the $L^p$-Sobolev spaces
$$
W^{s,p}(M,\Lambda^{q,r}):=  \{ u\in L^p(M,\Lambda^{q,r})\mid
\partial^\alpha u\in L^p\mbox{  for  }|\alpha |\le s\} ,
$$
for $1\le p \le \infty$, $s\in\NN$, where again differentiation is understood componentwise, with respect to local geodesic coordinates, and in the distributional sense.

As mentioned above, the group invariance and the compactness of the quotient
provide us with a number of useful uniformities. This applies, {\it e.g.}\ to
the pseudolocal estimates required in assumption (A) from \ref{assumption} above
in that all we will ever need will be derivable from the estimate for a single
neighborhood $U$ and a fixed pair of cutoffs $\zeta, \zeta'$, yielding a
universal $\epsilon>0$ and constant $C_{\zeta,\zeta'}$, as in \cite{PS}. We
refer the reader to \cite{D'A,C,E} for a discussion of this type of estimates as
well as sufficient geometric properties.

Let us end this section with a final word on forms and forms: Unfortunately we
need to  use these completely different concepts that bear the same name in
this paper. From Hilbert space theory we need \emph{sesquilinear forms} that are
\emph{bounded below}, e.g., the $Q=Q^{p,q}$ above. See Kato's \cite{K} and Reed
and Simon's \cite{RS}  classics and Faris' excellent lecture notes \cite{F} for
background. These forms are defined on $L^2$-spaces of \emph{differential
forms}, as we already mentioned. The standard reference for the relevant
notions of differential forms related to the $\dbar$-Neumann problem is
\cite{FK}.

%%%%%%%%%%%%%%%%%%%%%%%%%%%
\section{The intrinsic metric}
%%%%%%%%%%%%%%%%%%%%%%%%%%%
\label{metric_sec}
In \cite{PS} we used the intrinsic metric to bound the heat kernel of the
$\dbar$-Neumann Laplacian. Here it will again turn out to be extremely useful.
In this section we give a characterization and prepare the ground for a cutoff
procedure that is well suited to forms in the domain of the form $Q=Q^{p,q}$.
We rely on assumption \textbf{(A)} from \ref{assumption}, as usual.

\begin{definition}\label{metric} We define the $G$-invariant pseudo-metric on
$M$ by
\[ \dint (x,y) = \sup \{ w(y)- w(x) \mid w\in L^\infty\cap
C^\infty(\mbar,\mathbb R),\langle\dbar w,\dbar w\rangle_{\Lambda^{0,1}}\le
1\}.\]
\noindent
We define the distance between sets accordingly,
\[\dint (A;B):= \sup \{\inf_B w -\sup_Aw \mid w\in L^\infty\cap
C^\infty(\mbar,\mathbb R),\langle\dbar w,\dbar w\rangle_{\Lambda^{0,1}}\le 1\}\]
\noindent
for arbitrary $A,B\subset {\mbar}$. \end{definition}

Compared to the above definition, we extend the family of functions over which we take the supremum as follows.

\begin{lemma} Let $\mathcal A^1 =\{w\in C(\mbar,\mathbb R)\mid |\dbar
w|_{\Lambda^{0,1}}\le1,\, \mu-a.e.\}$ where the derivative is understood in the
distributional sense. It follows that any $w\in\mathcal A^1$ is a limit, locally
uniformly, of smooth functions $w_k$ with $|\dbar w_k|_{\Lambda^{0,1}}\le 1$.
\end{lemma}
\begin{proof} Apply Friedrichs mollifiers.
\end{proof}
\begin{corollary} $d_\square(x,y) = \sup\{w(y)-w(x)\mid w\in\mathcal A^1\}$.\end{corollary}
\begin{definition}\label{rhoE} For $E\subset\mbar$, put
\[\rho_E(x) =\inf\{d_\square(x,y)\mid y\in E\}.\]
\end{definition}
\begin{lemma} The function $\rho_E\in\mathcal A^1$ and $d_\square(\{x\}, E) =
\rho_E(x)$.\end{lemma}
We deduce the preceding lemma from the following description of the
saturation properties of $\mathcal A^1$.
\begin{prop}\cite[Props.\ A.1, A.2]{BdMLS} For $\mathcal A^1$ as above, we have the following properties.
\begin{enumerate}
\item $\mathcal A^1$ is \emph{balanced}, {\it i.e.}\ it is convex and closed under multiplication by -1.
\item $\mathcal A^1$ is closed under the operations $\min$ and $\max$
\item $\mathcal A^1$ is closed under the operation of adding constants
\item $\mathcal A^1$ is closed under pointwise convergence of functions uniformly bounded on compacts
\item Let $\mathcal F\subset\mathcal A^1\cap C(\mbar,\mathbb R)$ be stable under
$\max$ (resp.\ $\min$). If $u=\sup\{v\mid v\in\mathcal F\}$, (resp.\
$u=\inf\{v\mid v\in\mathcal F\}$) then $u\in\mathcal A^1$.
\end{enumerate}\end{prop}
\begin{proof} Note that the form $\mathcal E(u,v) = \langle\dbar u,\dbar v\rangle_{\Lambda^{0,1}}$ from \cite{PS} on $\mathcal D = \dom(Q^{0,0})$ is a strongly local Dirichlet form with energy measure
\[\Gamma(u,v) =  \langle\dbar u(x),\dbar v(x)\rangle_{\Lambda^{0,1}}\, d\mu(x),\]
\noindent
so the formalism of \cite[Appendix]{BdMLS} applies. \end{proof}

Now let us turn to an alternative description of the intrinsic metric. As calculated in \cite[\S3.2]{PS}, if the Hermitian metric on $\Lambda^{0,1}$ is associated to the Riemannian metric on $\Lambda^1$, then
\[ 2 \langle\dbar w,\dbar w\rangle_{\Lambda^{0,1}} =  \langle
dw,dw\rangle_{\Lambda^1}, \quad w\in C^\infty(\mbar,\mathbb R).\]
\noindent
By \cite{SC, Stu}, the form on the right induces the
Laplace-Beltrami operator $\Delta_{LB}$ on functions.

\begin{definition}\label{row} Where
\[L(\gamma) = \int_a^b |\dot\gamma(t)|_{T^1}\, dt\]
\noindent
for $\gamma$ piecewise smooth curves $\gamma:[a,b]\to\mbar$, and the length of
$\dot\gamma$ is measured with the Riemannian metric on $T^1M$, let
\[\rho(x,y) =  \inf \{L(\gamma)\mid \gamma\ {\rm is\ a\ piecewise\ smooth\ curve\ joining}\ x\ {\rm and}\ y\}.\]
 \end{definition}
\begin{corollary} In the situation above, we have
\begin{align} d_{\square} (x,y) =& \sqrt 2\sup\{w(y)-w(x) \mid  \langle
dw,dw\rangle_{\Lambda^1} \le 1\}\notag\\
=& \sqrt 2 \rho(x,y).\notag\end{align}
\end{corollary}
\begin{proof} Fix a $w\in\mathcal A^1\cap C^\infty(\mbar,\mathbb R)$ and let
$\gamma:[a,b]\to\mbar$ be a curve with
\[\rho(x,y)\ge L(\gamma) + \epsilon.\]
\noindent
We get
\begin{align} w(y)-w(x) =& w(\gamma(1))-w(\gamma(0)) = \int_0^1
\frac{d}{dt}w\circ\gamma(t)dt\notag\\ = &\int_0^1
\langle dw(\gamma(t)),\dot\gamma(t)\rangle dt  \notag\\
\le & \int_0^1 |dw(\gamma(t))|_{\Lambda^1}\,
|\dot\gamma(t)|_{T^1} dt \le  \sqrt 2 \int_0^1 |\dot\gamma(t)|_{T^1}
dt  \notag\\
\le & \sqrt 2 (\rho(x,y) + \epsilon).\notag \end{align}
\noindent
Thus $d_{\square} (x,y)\le\sqrt 2\rho(x,y)$. To show the reverse inequality, fix $y\in\mbar$ and note that it is
enough to prove that for $w(x) = \rho(x,y)$ we have weak differentiability and
 $|dw(x)|_{\Lambda^1}\le 1$, $\mu$-almost everywhere in $M$. By Rademacher's
theorem,
this amounts to showing that
\[|w(x) - w(x')| \le \rho(x,x')\]
\noindent
which follows from the triangle inequality for $\rho$.
 \end{proof}
\begin{rem} The existence of minimizing geodesics in the case at hand is demonstrated in \cite{AlexanderA-81}.
In the general Dirichlet form setting, the intrinsic metric gives at least a length space, as shown in \cite{Sto2}. From now on we will simply write $d(\cdot,\cdot)$ instead of $d_{\square}(\cdot,\cdot)$. \end{rem}
We now use the intrinsic metric to define cutoff functions.
Let $b>0$ and $\zeta\in C^1(\mathbb R, [0,1])$ so that
$\zeta|_{(-\infty,0]}\equiv 0$, $\zeta|_{[b,\infty)}\equiv 1$, and $\sup
|\zeta'(t)|<2/b$. It follows that
\begin{itemize}
 \item $\zeta\circ\rho_E\in W^{1,\infty}(\mbar,\mathbb R)$,
 \item $|\dbar(\zeta\circ\rho_E)|\le 2/b$,
 \item $\zeta\circ\rho_E|_{\bar E} \equiv
1$,
 \item $\zeta\circ\rho_E|_{B_b(E)^c} \equiv 0$.
\end{itemize}
\noindent
where $B_b(E) = \{y\in\mbar\mid d(y, E)\le b\}$ is the $b$-neighborhood of $E$.

A word on notation: For two quantities $A$ and $B$, we write $A\lesssim B$ to mean that there exists a constant $C>0$ such that $|A(\phi)|\le C|B(\phi)|$ uniformly for $\phi$ in whatever set relevant to the context.

We have the following elementary
\begin{lemma}\label{hat} For $u\in L^\infty (\mbar,\Lambda^k)$ with support in
$E$ and $v\in L^2(M, \Lambda^l)$, we have
\[\|u\wedge v\|_{L^2(M,\Lambda^{k+l})}^2 \lesssim \|u\|_{L^\infty(M,\Lambda^k)}^2 \int_E |v(x)|_{\Lambda^l}^2. \]
 \end{lemma}
\begin{proof} Pointwise, we have $|u\wedge v|_{\Lambda^{k+l}} \lesssim |u|_{\Lambda^k} |v|_{\Lambda^l}$, and by Sect.\ \ref{Prel}, $\|u\|_{L^\infty} = \esssup_M
|u|_{\Lambda^k}$, from which the result follows on integration.\end{proof}
\begin{prop} Let $\phi\in W^{1,\infty}(\mbar,\mathbb R)$ and $u\in\dom Q$. Then
$\phi u\in\dom Q$ and
\[Q(\phi u) \lesssim \|\phi\|_{W^{1,\infty}}^2 \left[ Q(u) + \|u\|_{L^2}^2\right].\]
\end{prop}
\begin{proof} First assume that $\phi$ and $u$ are smooth. Then
\[\dbar(\phi u) = \phi\dbar u+\dbar\phi\wedge u, \quad \dbar^*(\phi u) =  \phi\dbar^* u - \star[\partial\phi\wedge\star u],\]
\noindent
{\it cf.}\ \cite[Lemma 3.10]{PS}. The fact that the Hodge $\star$ is an isometry gives
\begin{align}Q(\phi u) =& \|\dbar(\phi u)\|_{L^2}^2 + \|\dbar^*(\phi u)\|_{L^2}^2\notag \\
 =&   \|\phi\dbar u\|_{L^2}^2 + \|\phi\dbar^*u\|_{L^2}^2 + \|\dbar\phi\wedge u\|_{L^2}^2 + \|\partial\phi\wedge\star u\|_{L^2}^2\notag \\
 & + 2\, \mathfrak{Re} \langle\phi \dbar u, \dbar\phi\wedge u\rangle_{L^2} - 2\, \mathfrak{Re} \langle\phi \dbar^* u, \star[\partial\phi\wedge \star u]\rangle_{L^2}.\notag
\end{align}
\noindent
With the previous lemma, Cauchy-Schwarz, and again the fact that $\star$ is an isometry, we obtain
\begin{align}\dots \le &\, \|\phi\|_{L^\infty}^2 Q(u) +  \|\dbar\phi\|_{L^\infty}^2 \|u\|_{L^2}^2 +  \|\partial\phi\|_{L^\infty}^2 \|u\|_{L^2}^2\notag\\
\notag & + \|\phi\dbar u\|_{L^2}^2 + \|\dbar\phi\wedge u\|_{L^2}^2 + \|\phi\dbar^* u\|_{L^2}^2+ \|\partial\phi\wedge \star u\|_{L^2}^2, \end{align}
\noindent
and each of these terms bounded by a constant multiple of the right-hand side in
the assertion. Now drop the assumption of smoothness and choose
$(\phi_k)_k\subset W^{1,\infty}\cap C^1$ so that $\phi_k\to\phi$, pointwise
a.e.\ and with $\|\phi_k\|_{W^{1,\infty}}\le\|\phi\|_{W^{1,\infty}}$, and
similarly $(u_k)_k\subset C_c^\infty(\mbar,\Lambda^{p,q})$ so that
\[Q(u_k-u)\longrightarrow 0 \quad {\rm and}\quad \|u_k-u\|_{L^2}\longrightarrow 0.\]
\noindent
It follows that $\phi_k u_k\to \phi u$ in $L^2$ and $\sup Q(\phi_k u_k)
<\infty$. Standard Fatou-type arguments \cite{MaR} give that $\phi u\in\dom Q$
and
\[Q(\phi u)\le \lim\inf Q(\phi_k u_k), \]
\noindent
giving the assertion.\end{proof}

%%%%%%%%%%%%%%%%%%%%%%%%%%%%%%%%
\section{The Caccioppoli inequality}
%%%%%%%%%%%%%%%%%%%%%%%%%%%%%%%%
\label{Sec_4} As usual, we work under the assumption \textbf{(A)} from
\ref{assumption} above. Let us first introduce the notion of a generalized
eigenform.
\begin{definition} A form $u\in L^2_{\rm loc}(M,\Lambda^{p,q})$ is said to be a
\emph{generalized eigenform} for $\square = \square_{p,q}$  if

\medskip

\begin{flushleft}
1) $u\in{\rm dom}_{\rm loc}\, Q$. That is, for any compact $K\subset\mbar$
there is a $v\in\dom Q$ such that $v|_K = u|_K$ and

\medskip

2) There exists a $\lambda\in\mathbb R$ such that $Q(u,\phi) = \lambda\langle
u,\phi\rangle$ for all $\phi\in C_c^\infty(\mbar,\Lambda^{p,q})$.
\end{flushleft}
\end{definition}

\begin{rem}  Note that ${\rm dom}_{\rm loc} Q$ is in $H^{1/2}_{\rm loc}$; see \cite[Prop.\ 1.2]{GHS}. Note also that the identity in 2) is a weak form of the equation $\square u=\lambda u$.

By locality of the energy we can define $Q(u,\phi)=Q(v,\phi)$
provided $\supp\phi\subset K$ and $v$ is as in the definition. Alternatively, we
can write
\[Q(u,\phi) = \int \langle\dbar u(x),\dbar \phi(x)\rangle_{\Lambda^{p,q+1}} +
\langle\dbar^* u(x),\dbar^* \phi(x)\rangle_{\Lambda^{p,q-1}} \ d\mu(x)\]
\noindent
noting that the integral is convergent. Moreover, we have that
$$
u\in{\rm dom}_{\rm loc}\, Q \Leftrightarrow \dbar u \in L^2_{\rm loc}(M,
\Lambda^{p,q+1})\quad \mbox{  and  } \quad \dbar^* u \in L^2_{\rm loc}(M,
\Lambda^{p,q-1}).
$$
\end{rem}
The Caccioppoli inequality states that for any  generalized
eigenform $u$, the
energy
\[M\ni x\longmapsto \langle\dbar u(x),\dbar u(x)\rangle_{\Lambda^{p,q+1}} +
\langle\dbar^*
u(x),\dbar^* u(x)\rangle_{\Lambda^{p,q-1}}\]
\noindent
is locally bounded by the $L^2$-norm of $u$.
We follow the strategy of \cite{BdMLS} in what follows. See also \cite[Prop.\ 3]{AV}, for a similar result. The authors of the former
paper had not been aware of the latter at the time their paper appeared. The
reader should note one important difference between the result in  \cite{AV} and
in the following result: $u$ is not supposed to be in the
domain of the operator, or even locally in the Sobolev space $H^2$.

\begin{theorem}{\bf (Caccioppoli inequality)} For any generalized eigenform
$u$ of $\square_{p,q}$ associated to an eigenvalue $\lambda\ge 0$, every
compact set $E\subset\mbar$, and every $b\in(0,1]$ we have
\[\int_E |\dbar u|_{\Lambda^{p,q+1}}^2 + |\dbar^*
u|_{\Lambda^{p,q-1}}^2\le
2\lambda\int_{E}
|u|_{\Lambda^{p,q}}^2 +
\frac{4}{b^2}\int_{B_b(E)}
|u|_{\Lambda^{p,q}}^2. \]
\end{theorem}
\begin{rem}Note that in order to control the energy on $E$ we need to take the
$L^2$-norm on a slightly larger set $B_b(E)$, the $b$-neighborhood of $E$.
\end{rem}
\begin{proof} Pick a cutoff function $\eta = \zeta\circ\rho_E$ as constructed in
Sect.\ \ref{metric}, so that $|\dbar\eta|\le 2/b$, $\eta|_{\bar E} \equiv 1$,
and $\eta|_{B_b(E)^c} \equiv 0$.
\noindent
The eigenvalue equation
\[Q(u,\phi) = \lambda\langle u,\phi\rangle_{L^2(M,\Lambda^{p,q})}, \qquad
(\phi\in C_c^\infty(\mbar,\Lambda^{p,q})),\]
\noindent
extends to arbitrary $\phi\in\dom Q$ by approximation. Therefore we may
calculate
\begin{align} \lambda\langle u,\eta^2 u\rangle = & \, Q(u,\eta^2 u)\notag \\
= & \int_{B_b(E)} \langle\dbar u, \dbar(\eta^2 u)\rangle + \langle\dbar^* u,
\dbar^*(\eta^2 u)\rangle \notag\\
= & \int_{B_b(E)} \eta^2 [|\dbar u|^2 + |\dbar^* u|^2] + \langle\dbar
u,\dbar\eta^2\wedge u\rangle - \langle\dbar^* u, \star[\partial\eta^2\wedge\star
u]\rangle. \notag \end{align}
\noindent
Leibniz' rule gives
\[\dots = \int_{B_b(E)} \eta^2 [|\dbar u|^2 + |\dbar^* u|^2] + 2\langle\eta\dbar
u,\dbar\eta\wedge u\rangle - 2\langle\eta\dbar^* u,
\star[\partial\eta\wedge\star u]\rangle. \]
\noindent
Now rearrange terms, apply Cauchy--Schwarz, and Lemma \ref{hat} to get
\begin{align} \int_{B_b(E)} \eta^2 [|\dbar u|^2 + |\dbar^* u|^2]\le &\,
\lambda\langle u,\eta^2 u\rangle + 2\|\eta\dbar u\|\|\dbar\eta\wedge u\| + 2
\|\eta\dbar^* u\|\|\partial\eta\wedge\star u\| \notag\\
\le & \, \lambda \int_{B_b(E)} \eta^2 |u|^2  + \frac12 \left[\|\eta\dbar^* u\|^2
+
4 \|\partial\eta\wedge\star u\|^2\right]. \notag \end{align}
The second term on the right hand side is 1/2 the left hand side so we have
\[\frac12  \int_{B_b(E)} \eta^2 [|\dbar u|^2 + |\dbar^* u |^2]\lesssim
\lambda\int_{B_b(E)} |u(x)|_{\Lambda^{p,q}}^2 + \frac{2}{b}
\int_{B_b(E)\setminus E}  |u(x)|_{\Lambda^{p,q}}^2, \]
\noindent
which yields the assertion since $\eta\equiv 1$ on $E$.\end{proof}

%%%%%%%%%%%%%%%%%%%%%%%%%%%%%%%%
\section{Subexponentially bounded eigenforms induce spectrum}
%%%%%%%%%%%%%%%%%%%%%%%%%%%%%%%%
\label{Sec_5}
Here we closely follow \cite{BdMLS}, see also the survey \cite{LSV1}. The
treatment here rests on two main observations. The first is a criterion for
$\lambda\in\sigma(H)$ in terms of the quadratic form $h$ associated with $H$.
\emph{Singular sequences} or \emph{Weyl sequences} for $H$ and $\lambda$ are
sequences $(f_n)_{n\in\NN}\subset \dom(H)$ that satisfy $\| f_n\|=1$ for all
$n\in\NN$ and
$$
\| Hf_n-\lambda f_n\|\to 0\mbox{  for  }n\to\infty .
$$
Clearly, the existence of such a sequence implies that $\lambda\in\sigma (H)$.
However, due to the requirement $f_n\in \dom(H)$ such singular sequences may be
hard to construct. The next proposition gives a quadratic form version, which
clearly is easier to find. Note that the terms ``singular sequence'' and ``Weyl
sequence'' are most commonly used in a stricter sense; namely, it is required
that, additionally, $f_n\stackrel{w}{\to}0$. In this case one even gets that
$\lambda$ lies in the essential spectrum of $H$. Our criterion below does not
require this. E.g., for an eigenvalue $\lambda$ one could simply take $f_n=f$,
where $f$ is a normalized eigenelement of $H$ with eigenvalue $\lambda$.

 Throughout this section we assume (A) from
\ref{assumption}.

\begin{prop}{\bf (Weyl type Criterion)} Let $h$ be a closed semibounded
form and let $H$ be the associated self-adjoint operator. Then the following are
equivalent:

\medskip

\noindent
1) $\lambda\in\sigma(H)$

\medskip

\noindent
2) There exists a sequence $(u_k)_k$ in $\dom h$ with $\|u_k\|\to 1$ and
\[\sup \{ \left|h(u_k,v)-\lambda\langle u_k,v\rangle_{L^2}\right| \mid v\in\dom
h,\, \|v\|_h \le 1 \}
\longrightarrow 0 \]
\noindent
for $k\to\infty$. \end{prop}
For the proof see \cite[Lemma 1.4.4]{Sto2} and \cite{DDI}.

As a last ingredient for the main result, let us introduce the {\it inner
b-collar} of a set $E\subset M$, given by
\[\mathcal C_b(E) = \{x\in E\mid d(x,E^c)\le b\}.\]
\begin{theorem}{\bf (1/2-Schnol)}\label{half} Assume that $\lambda\in\mathbb R$
admits a
generalized eigenform $u$ so that there exists a sequence $E_k$ of compact
subsets of $\mbar$ and $b>0$ with
\[\frac{\|u \bid_{C_b(E_k)}\|}{\|u \bid_{E_k}\|} \longrightarrow 0 \quad {\rm
as} \quad k\to\infty.\]
\noindent
Then $\lambda\in\sigma(H)$.
\end{theorem}
\begin{proof} Let us first calculate, for $\eta\in W_c^{1,\infty}(\mbar,\mathbb
R)$ and $u,v\in{\rm dom}_{\rm loc}\, Q^{p,q}$,
\begin{align}Q(\eta u,v) & -  Q(u,\eta v)\notag\\
=&\,  \int_M \langle\dbar(\eta u),\dbar v\rangle -\langle\dbar
u,\dbar(\eta v)\rangle + \langle\dbar^*(\eta u),\dbar^* v\rangle -\langle\dbar^*
u,\dbar^*(\eta v)\rangle \label{qdiffs} \\
 =&  \int_M  \langle\dbar\eta\wedge u,\dbar v\rangle -\langle\dbar
u,\dbar\eta\wedge v\rangle +\dots \notag\\
 &\qquad +  \langle\star[\partial\eta\wedge\star u],\dbar^* v\rangle
-\langle\dbar^* u,\star[\partial\eta\wedge\star v]\rangle.\notag\end{align}
Now choose a sequence $E_k$ as in the assumptions and define
\[F_k = \{x\in E_k\mid d(x,E_k^c)\ge b/2\},\]
\noindent
with which we will define suitable cutoff functions. So pick $\zeta\in
C^1(\mathbb R)$ with $0\le\zeta\le 1$, $\zeta|_{(-\infty,0]}\equiv 1$,
$\zeta|_{[b/4,\infty)}\equiv 0$, and $\sup|\zeta'|\le 8/b$. Note that
$\eta_k:=\zeta\circ\rho_{F_k}\in W_c^{1,\infty}(\mbar,\mathbb R)$ satisfies
$0\le\eta_k\le\bid_{B_{b/4}(F_k)}$ and $\supp|\partial\eta_k|\subset
B_{b/4}(F_k)\setminus F_k =:G_k$. Moreover, note that
$B_{b/4}(G_k)\subset
C_b(E_k)$.

We now show that
\[u_k = \frac{\eta_k u}{\|\eta_k u\|_{L^2}}\]
\noindent
gives an approximate eigensequence as required by the Weyl criterion above.

Let $v_k=\eta_k u$. For $v\in\dom Q^{p,q}$ with $\|v\|_{Q}\le 1$, we estimate
\begin{align}Q(u_k,v) - \lambda\langle u_k,v\rangle =&\,
\frac{1}{\|v_k\|_{L^2}}\left[Q(v_k,v) - \lambda\langle
v_k,v\rangle\right]\notag\\
=&\, \frac{1}{\|v_k\|_{L^2}}\left[Q(\eta_k u,v) - \lambda\langle u,\eta_k
v\rangle\right]\notag\\ =&\, \frac{1}{\|v_k\|_{L^2}}\left[Q(\eta_k u,v) -
Q(u,\eta_k v)\right],\notag\end{align}
\noindent
since $\eta_k$ is real-valued. We have used that $u$ is a generalized eigenform
with eigenvalue $\lambda$ and the fact, discussed in the proof of Caccioppoli's
inequality, that $\eta_k v$ can be taken to be a test function. Now, as in
\eqref{qdiffs},
\begin{align}Q(u_k,v) & - \lambda\langle u_k,v\rangle\notag\\
 =&\, \frac{1}{\|v_k\|_{L^2}}\int_M\langle\dbar\eta_k\wedge u,\dbar
v\rangle - \langle\dbar u,\dbar\eta_k\wedge  v\rangle \notag\\
 &\quad-\langle\star[\partial\eta_k\wedge\star u],\dbar^* v\rangle +
\langle\dbar^* u,\star[\partial\eta_k\wedge\star v]\rangle. \notag\end{align}
\noindent
Due to the support properties of the $\eta_k$, we know that
\begin{align}\dots& \lesssim \frac{1}{\|v_k\|_{L^2}}\int_{G_k}
\|\dbar\eta_k\|_\infty \left[\|u\|_{L^2}\|\dbar
v\|_{L^2}+\|u\|_{L^2}\|\dbar^*v\|_{L^2} \right] \notag\\
&\qquad +  \frac{1}{\|v_k\|_{L^2}}\int_{G_k} \|\dbar\eta_k\|_\infty
\left[\|\dbar u\|_{L^2}\| v\|_{L^2}+\|\dbar^*u\|_{L^2}\|v\|_{L^2} \right]
\notag\\
&\lesssim  \frac{1}{\|v_k\|_{L^2}} \left[\|u \bid_{G_k}\|_{L^2} \|v\|_Q +
\|\dbar u\|_{L^2(G_k)} \|v\|_{L^2(M)} +
\|\dbar^*u\|_{L^2(G_k)}\|v\|_{L^2(M)}\right].\notag\end{align}
\noindent
Now apply Caccioppoli with $E=G_k$ and $b/4$ to get
\begin{align}\dots &\lesssim  \frac{1}{\|v_k\|_{L^2}} \left[\|u
\bid_{G_k}\|_{L^2} + \| u\|_{L^2(B_{b/4}(G_k))} \right]\notag\\
&\lesssim
\frac{\|u\bid_{C_b(E_k)}\|_{L^2}}{\|u\bid_{E_k}\|_{L^2}}\notag\end{align}
\noindent
since $0\le \eta_k\le \bid_{E_k}$ and $B_{b/4}(G_k)\subset C_b(E_k)$.
\end{proof}
Generalizing the notion from statistical mechanics, let us call a sequence
$(E_k)$ a \emph{van Hove sequence} if it has the property that
$$\frac{\vol C_b(E_k)}{\vol E_k}\longrightarrow 0\qquad {\rm as} \qquad
k\longrightarrow \infty,
$$
for some $b>0$.
\begin{corollary} Assume that $\mbar$ admits a van Hove sequence. It follows
that
$0\in\sigma(\square_{0,0})$ and $1$ is a generalized eigenfunction.\end{corollary}
\begin{proof}Clearly, $1$ is a generalized eigenfunction for the eigenvalue $0$.
By
the preceding remark, it satisfies the requirement for the theorem
above.\end{proof}
Apart from certain uniformities, the $G$-invariance which we assume
throughout this treatment is certainly too strong a condition to impose.  Note
that a suitable generalization of the theorem above allows for manifolds with
very different geometries in different ``directions to infinity.''  One such
direction which supports a van Hove sequence  is sufficient for 0 to be in the
spectrum of $\square_{0,0}$.

We now add some sufficient conditions for the assumptions in the theorem. They
rest on the following notions: A function $J:[0,\infty)\to[0,\infty)$ is said to
be {\it subexponentially bounded} if for any $\alpha>0$ there exists a
$C_\alpha>0$ such that
\[J(r)\le C_\alpha e^{\alpha r} \qquad (r\ge 0).\]
Similarly, a form $u\in L^2_{\rm loc}(M,\Lambda^{0,1})$ will also be called
subexponentially bounded if for some $z_0\in\mbar$,
\[e^{-\alpha w}u\in L^2(M,\Lambda^{p,q})\]
\noindent
for any $\alpha>0$, where $w(z) = d(z,z_0)$.

As in lemmata 4.2, 4.3, and Thm.\ 4.4 of \cite{BdMLS}, we obtain
\begin{corollary} \label{sub}
Assume that $\lambda\in\mathbb R$ admits a
subexponentially
bounded eigenform for $\square$. It follows that
$\lambda\in\sigma(\square)$.\end{corollary}
This or the previous corollary has as a special case
\begin{corollary} Assume that there is a $z_0\in\mbar$ such that $r\mapsto\vol
B_r(z_0)$ is subexponentially bounded.  It follows that
$0\in\sigma(\square)$.\end{corollary}
\begin{rem} See the example in \cite{DSP}.\end{rem}
During the writing of \cite{BdMLS,LSV1} we were not aware of M.\ Shubin's papers
\cite{Shu1,Shu2}, where strongly related results are presented. The main
difference is that our approach is based on the underlying forms, making it
applicable in cases where nothing is known about the domain of the operator. On
the other hand, the latter papers contain results about higher order elliptic
operators.

%%%%%%%%%%%%%%%%%
\section{Expansion in generalized eigenforms}
%%%%%%%%%%%%%%%%%%%
\label{Sec_6}
Here we prove Theorem \ref{thm1.3}, in fact the stronger result Proposition
\ref{eigen-prop} below, where assumption (A) from \ref{assumption} is required,
as usual.

Some explanations are in order: {\it spectrally a.e.}\ means a.e.\ with respect
to a spectral measure; in turn, a spectral measure $\rho$ is a measure with the
property that $\rho (I)=0$ if and only if
$E_I(\square)= 0$, where $E_\cdot(\square)$ denotes the spectral projection
of the operator $\square$.

The strategy of proof is sufficiently parallel to the one in \cite{BdMS} so that
we do not carry out all the details but rather point at differences; we fixed
integers $p\ge 0, q\ge 0$ so that the pseudolocal estimate holds true. This
latter condition is important in that we use
ultracontractivity established in \cite{PS}, i.e., $e^{-t\square}: \lzpq\to \lupq$ for $t>0$. The
compactness property referred to above is contained in the following:

\begin{lemma} In the situation of the theorem above let $\gamma(x):=e^{-tx}$ and
$T:=M_{\omega^{-1}}$ the multiplication operator. Then $\gamma(\square)T^{-1}$
is Hilbert-Schmidt. \end{lemma}
\begin{proof} This follows from the factorization principle based on
Grothendieck's theorem. See \cite{DF} for the abstract background and \cite{BdMS}
for an application in a situation similar to ours.

Indeed, for bounded operators, from
\[A:L^2 \longrightarrow L^\infty, \quad B: L^\infty\longrightarrow L^2, \]
\noindent
it follows that $BA:L^2\to L^2$ is a Hilbert-Schmidt operator. We can apply this to deduce that
\[(\gamma(\square)T^{-1})^* = (T^{-1})^* \gamma(\square)^*\]
\noindent
is Hilbert-Schmidt: $\gamma(\square):L^2\to L^\infty$ is the above mentioned uultracontractivity and $T^{-1} = M_\omega: L^\infty \to L^2$, since
$\omega$ is an $L^2$ function. Since the adjoint of a Hilbert-Schmidt operator is likewise
Hilbert-Schmidt, we have the result. \end{proof}
Suppressing the indices $p,q$, let
\[\cH_+:=\{ \alpha\in \lzpq\mid \alpha\in \dom(T)\},\]
\noindent
and $\cH_-$, the completion of $\cH:=\lzpq$ with respect to the inner product
$\langle \alpha ,\beta \rangle_-:= \langle
T^{-1}\alpha,T^{-1}\beta\rangle_{\cH}$.
We have a special case of a Gelfand triple here, considering on $\cH_+$ the
inner product $\langle \alpha ,\beta \rangle_+:= \langle
T\alpha,T\beta\rangle_{\cH}$.
\begin{rem} We have that
$$
C^\infty_c(M,\lpq)\subset
\{ \alpha\in \dom (\square)\cap\dom(T)\mid \square\alpha\in\dom(T)\}
$$
is dense in $\cH$. In the next section we will prove much more, namely
that $\ccpq$ is a core for $\square$. Note the important difference between $M$ and $\mbar$ here.
\end{rem}
%%%%%%%

In the following result we see a much stronger though more technical version of
the theorem above. It uses the notion of an \emph{ordered spectral
representation}, that goes as follows: Given is a self adjoint operator $H$ in
some Hilbert space $\cH$, a spectral measure $\rho$ of $H$,
$N\in\NN\cup\{\infty\}$ a sequence $(M_j)_{j<\infty}$ of measurable subsets
$M_j\subset \RR$ so that $M_j\supset M_{j+1}$ and a unitary
\[U=(U_j)_{j<\infty}:\cH\to\bigoplus_{j< N} L^2(M_j,\rho)L^2(M_j,d\rho)\]
\noindent
so that
\[U\varphi(\square)=M_\varphi U\]
\noindent
for every bounded measurable function $\varphi$ on $\RR$.
\begin{prop} \label{eigen-prop}
Let $\rho$ be a spectral measure for $\square$
and $U=(U(j))_{j<
N}$, $N\in\NN\cup\{\infty\}$, an ordered spectral representation for $\square$.
Also let $\omega$, $T$, $\cH_+$ and $\cH_-$ be as above. Then there are
measurable functions $M_j\to\cH_-,\lambda\mapsto\varepsilon_{j,\lambda}$ for
$j\in\NN$,
$j<N$ such that:
\begin{itemize}
 \item[{\rm (1)}]
$U_j\alpha(\lambda)=\langle\alpha,\varepsilon_{j,\lambda}\rangle$ for
$\alpha\in\cH_+$ and $\rho$-a.e.\ $\lambda\in M_j$.
 \item[{\rm (2)}]
For every $g=(g_j)_{j< N}\in\bigoplus_{j< N} L^2(M_j,\rho)$ we have
\[U^{-1}g=\lim_{m\to N, R\to\infty}\sum_{j=1}^m\int_{M_j\cap
[-R,R]}g_j(\lambda)\varepsilon_{j,\lambda}d\rho(\lambda),\]
and therefore, for every $\alpha\in\cH$,
\[\alpha = \lim_{m\to N, R\to\infty}\sum_{j=1}^m\int_{M_j\cap
[-R,R]}U_j\alpha(\lambda)\varepsilon_{j,\lambda}d\rho(\lambda).\]
 \item[{\rm (3)}] If $\alpha\in \dom(\square)\cap\cH_+$ with
$\square\alpha\in\cH_+$, then
\[\langle\square\alpha,\varepsilon_{j,\lambda}\rangle
=\lambda\langle\alpha,\varepsilon_{j,\lambda}\rangle\mbox{  for }\rho-a.e.\
\lambda\in M_j.\]
\end{itemize}
\end{prop}

For details on ordered spectral representations, see \cite{PSW}; this reference
is the basis for our proof of the eigenform expansion.

Part 3 of the above proposition ensures that
\[\square \varepsilon_{j,\lambda} = \lambda \varepsilon_{j,\lambda}\]
\noindent
in the weak sense. This is why we speak of a {\it generalized eigenform}. 
\begin{rem} Due to the interior ellipticity of $\square$, \cite[Thm.\ 2.2.9]{FK} we obtain that the eigenforms constructed above are in $C^\infty(M,\Lambda^{p,q})$ for $q>0$.\end{rem}

%%%%%%%%%%%%%%%%%%%%%%%%%%%%%%%%

%%%%%%%%%%%%%%%%%%%%%%%%%%%%%%%%
\section{Essential self-adjointness of $\square$}\label{essential}
%%%%%%%%%%%%%%%%%%%%%%%%%%%%%%%%
As we explained in the introduction, $\square$ is defined via its sesquilinear
form, so its domain $\dom(\square)$ is only given implicitly. In the previous
sections we have seen that even without explicit knowledge of its domain we can
analyze important properties of $\square$.

On the other hand it is known for manifolds without boundary that
elliptic operators are typically essentially self-adjoint on smooth compactly
supported forms, see {\it e.g.}\ \cite{Shu1,Shu2,Shu3} and the literature
cited there. Thus it is a natural question
whether the same holds true in the situation at hand with two important
differences: there is a boundary, and we do not have ellipticity but only
subellipticity.

Essential self-adjointness means that there is a unique self-adjoint extension
of $\square|_{\domc}$ and this is in turn  equivalent  to the fact that
$\domc:=\domc\square:= \dom(\square)\cap \ccpq$ is a \emph{core for}
$\square$, i.e., $\overline{\square|_{\domc}}=\square$, where $\overline{T}$
denotes, as usual, the closure of the operator $T$. We want to point out that
there is a big difference due to  the boundary: in the usual complete case
without boundary, the so-called \emph{minimal operator}, defined on $\ccpq$ is
essentially self-adjoint. This fails in our situation. There are various
different self-adjoint extensions. E.g., the operator $\square_{p,q}$ we consider
is obviously different from the Friedrichs extension of
$\square_{p,q}|_{C^\infty_c(M,\lpq)}$ which would usually be called the
$\square$ with Dirichlet boundary conditions, and which has a smaller form
domain.

A first step in showing the asserted essential self-adjointness is the following
result from  \cite{GHS}. As usual, asumption \textbf{(A)} from
\ref{assumption} is in force.
Here and for what follows we fix $\rho$ to be the (positive) distance to $bM$ as given by a $G$-invariant Riemannian metric on $M$.

\begin{prop}\label{dbarlap} Let $\vartheta$
be the formal adjoint operator to $\dbar$, and denote by
$\sigma=\sigma(\vartheta,\cdot)$ its principal symbol. Assume also that
 $q>0$ and let $\square=\square_{p,q}$. Then
\begin{align}\domo\square := \{u\in C^{\infty}(M, \Lambda^{p,q}) \mid &
u,\dbar u, \vartheta u\in L^2,\notag\\ & \sigma(\vartheta,d\rho)u|_{bM}=0,\
\sigma(\vartheta,d\rho)\dbar u|_{bM}=0\}.\notag\end{align}
is a core for $\square$.
\end{prop}

\begin{proof} For convenience, we reprove this proposition here. Let
$u\in\dom\square_{p,q}$. Then $\square u + u = \alpha\in L^2$. Now let
$(\alpha_k)_k \subset C^\infty_c(\mbar, \Lambda^{p,q})$, so that $\alpha_k\to
\alpha$ in $L^2$ and put $u_k = (\square + \bid)^{-1}\alpha_k$. Since $(\square
+ \bid)^{-1}$ is defined everywhere and is bounded, we have that $(u_k)_k$ is
Cauchy in $L^2$ with limit $u$. Applying the pseudolocal estimate \cite{E},
\cite[Thm.\ 2.4]{PS} we have
\[\|\zeta u_k\|_{H^{s+\varepsilon}}=\|\zeta (\square
+\bid)^{-1}\alpha_k\|_{H^{s+\varepsilon}}\lesssim
\|\zeta'\alpha_k\|_{H^s}+\|\alpha_k\|_{L^2},\]
\noindent
thus $(u_k)_k\subset C^\infty(\mbar,\Lambda^{p,q})$ and we have shown that
the assertion is true.
\end{proof}

For the proof of essential self-adjointness we need some geometrical tools.
First recall the following
\begin{definition}\cite[p33]{FK}, \cite[\S2.2]{Stra} A \emph{special boundary chart} $U$ is a chart intersecting $bM$ having the following properties:
\begin{enumerate}
\item With $\rho$ the function defining $bM$ as above, the functions $t:=\{t_1,...,t_{2n-1}\}$, together with $\rho$ form a coordinate system on $U$.
\item The functions $\{t, \rho=0\}$ form a coordinate
system on $bM\cap U$.
\item With respect to the Riemannian structure in the cotangent bundle, choose a local orthonormal basis $\omega_1,\dots,\omega_n$ for $C^\infty(U,\Lambda^{1,0})$ such that $\omega_n=\sqrt{2}\ \partial\rho$ on $U$.
\end{enumerate}\end{definition}

 Let us describe $\dom\square$ by restating the boundary conditions as in
\cite[\S5.2]{FK}. In terms of the Hermitian structure $\langle\ , \,
\rangle_\Lambda$ in $\Lambda^{p,q}$, the above conditions on the symbol
$\sigma(\vartheta,d\rho)$ translate to the following criteria. Members of
$\domo\square$ are those forms $\phi\in C^\infty(\bar
M,\Lambda^{p,q})$ satisfying the following
$\dbar$-Neumann boundary conditions
\begin{enumerate}
\item $\langle\phi,\dbar\rho\wedge\psi\rangle_\Lambda |_{bM} = 0,\quad
(\psi\in\Lambda^{p,q-1})$,  and
\item $\langle\dbar\phi,\dbar\rho\wedge\psi\rangle_\Lambda |_{bM} = 0,\quad
(\psi\in\Lambda^{p,q})$.
\end{enumerate}
The first condition (equivalent to $\phi\in\domo\vartheta$) is obviously
preserved by introduction of a cutoff function $\phi\to\chi\phi$ since the condition is
algebraic.

The second ``free boundary'' condition becomes
\[ \langle\dbar(\chi\phi),\dbar\rho\wedge\psi\rangle_\Lambda |_{bM} =
\langle(\dbar\chi)\wedge\phi,\dbar\rho\wedge\psi\rangle_\Lambda |_{bM} +
\langle\chi\dbar\phi,\dbar\rho\wedge\psi\rangle_\Lambda |_{bM} = 0.\]
\noindent
Upon restriction to the boundary, the second term is zero by assumption that
$\phi\in\dom\square$, which assumes that $\dbar\phi\in\dom\dbar^*$. Thus we are interested in the condition
\[\langle(\dbar\chi)\wedge\phi,\dbar\rho\wedge\psi\rangle_\Lambda|_{bM} = 0,\
\forall\psi\in\Lambda^{p,q}.\]
In terms of the forms defined in the special boundary chart, we have the formulas
\[\dbar\rho = \frac{1}{\sqrt 2}\bar\omega^n, \qquad \dbar\chi = \bar L_k\chi\,
\bar\omega^k\]
\noindent
so cutoff functions $\chi$ satisfying
\begin{equation}\label{tazera}\bar L_n \chi|_{bM} = 0\end{equation}
\noindent
preserve $\dom\square$. Notice that there are no other restrictions on $\chi\in
C^\infty(\mbar)$ beyond
this one at the boundary, so $\chi$ satisfying \eqref{tazera} may be extended smoothly to the interior of $M$ in an arbitrary way.

We may write the relation \eqref{tazera} in such a way
that manifestly separates the tangential and normal derivatives of $\chi$, as indicated in \cite[p86]{BS}. First
note that since $L_n$ is dual to $\omega^n = \sqrt 2\partial\rho$, we have
\[L_n \rho = d\rho(L_n) = \langle
(\partial+\dbar)\rho,\omega^n\rangle_\Lambda  = \sqrt
2\langle\partial\rho,\partial\rho\rangle_\Lambda,\]
\noindent
and similarly $\bar L_n \rho = \sqrt
2\langle\dbar\rho,\dbar\rho\rangle_\Lambda$. It follows that $(\bar L_n -
L_n)\rho=0$ and thus $\bar L_n - L_n$ is a vector field tangential to $bM$. If $J$ is the complex structure, then $L_n$ and $\bar L_n$ lie in
the $i$ and $-i$ eigenspaces of $J$, respectively and
\[J(\bar L_n - L_n) = -i(\bar L_n + L_n)\]
\noindent
must not be tangential; indeed, $(\bar L_n + L_n)\rho = 2\sqrt
2\langle\dbar\rho,\dbar\rho\rangle\neq 0$. The same calculations provide that
the equation
\begin{equation}\label{cond}-iJ(\bar L_n - L_n)\chi = (\bar L_n - L_n)\chi\end{equation}
\noindent
(in $bM$) is equivalent to the property $\bar L_n\chi|_{bM}=0$. Since only the
normal derivative is prescribed at the boundary, it follows that given any
smooth function $\chi$ in $bM$, there exists an extension to a collar of $bM$
which fulfills the requirement \eqref{tazera}. Thus, {\it any} function $\chi\in
C^\infty(bM)$ can be extended to $\mbar$ in such a way that \eqref{tazera}
holds, {\it cf.}\ Lemma \ref{fundamental} below.

\begin{definition} A sequence of functions $(\chi_k)_k$ in $C^\infty_c(\mbar,
\RR)$ is called a \emph{good cutoff-exhaustion
of} $\mbar$ if
\begin{itemize}
 \item [(C1)] $\chi_k\to 1$ as $k\to\infty$,
 \item [(C2)] $\bar L_n \chi_k|_{bM} = 0$ for all $k\in\NN$,
and
 \item [(C3)] $\sup\{ \| \partial^\alpha\chi_k\|_\infty, |\alpha|\le
m\}<\infty$, for any $m\in \NN$,
 \end{itemize}
\end{definition}
\noindent
where the derivatives in the last condition are with respect to geodesic coordinates. 
Note that $\bar L_n$ is globally defined in a collar of the boundary of $M$.

Our goal here will be to demonstrate the existence of good cutoff-exhaustions of
$\mbar$ and to use such a sequence can to show that $\domc$ is a
core for $\square$. We start with:

\begin{prop}\label{prop7.4}
 Let $U$ be a special boundary chart and $\chi\in
C_c^{\infty}(U,\RR)$ with $\bar L_n \chi|_{bM} = 0$. Then, for any
$u\in\dom\square_{p,q}$ for $q>0$ it follows that
$$\chi u\in\dom\square$$
and
\begin{equation}\label{chi_u}
\|\square(\chi u)\|_{L^2}\lesssim\underbrace{\sup\{ \|
\partial^\alpha\chi\|_\infty, |\alpha|\le 2
\}}_{=\|\chi\|_{W^{2,\infty}}}\cdot (\|\square u\|^2_2 +\| u\|^2_2)^{\frac12}.
\end{equation}

\end{prop}
\begin{proof}
The factor $(\|\square u\|^2_2 +\| u\|^2_2)^{\frac12}$ appearing above is called the operator norm $\| u\|_{\square}$ of $u$. It dominates
the form norm $\| u\|_Q$ in the sense that
$$\| u\|_Q:=(Q(u,u) +\| u\|^2_2)^{\frac12}\lesssim \|
u\|_{\square} .$$
 Let us first consider the case that $u\in\domo\square_{p,q}$. By the
calculation above, $\chi u\in\dom\square$. For the proof of the estimate
(\ref{chi_u}), we
use the following straightforward calculation:
\begin{eqnarray*}
 \dbar(\chi u)&=& (\dbar \chi)\wedge u + \chi \dbar u \\
\dbar^\ast(\chi u)&=&(-\star \partial \star)(\chi u)\\
&=& -\star[(\partial\chi)\wedge\star u]+\chi \dbar^\ast u,
\end{eqnarray*}
from which we get that
$$
 |\langle \square(\chi u),v\rangle - \langle \chi\square u,v\rangle|=
|Q(\chi u,v)-Q(u,\chi v)|$$
$$
= | \langle \dbar(\chi u),\dbar v\rangle+ \langle \dbar^\ast(\chi
u),\dbar^\ast v\rangle - \langle \dbar u,\dbar(\chi v)\rangle
-\langle \dbar^\ast u,\dbar^\ast(\chi v)\rangle| $$
$$
= | \langle \dbar\chi\wedge u,\dbar v\rangle -\langle \dbar u,\dbar\chi\wedge
v\rangle -\langle \star(\partial\chi\wedge \star u),\dbar^\ast v\rangle +
\langle \dbar^\ast u,\star(\partial\chi\wedge \star v)\rangle | .
$$
The first term can be estimated as
$$
| \langle \dbar^\ast(\dbar\chi\wedge u), v\rangle|\lesssim
\|\chi\|_{W^{2,\infty}}\| u\|_Q \|v\|
$$
and similarly we can bound the third term. The second and fourth terms are
easily bounded and we get
$$
 |\langle \square(\chi u),v\rangle - \langle \chi\square u,v\rangle |
\lesssim
\|\chi\|_{W^{2,\infty}} \| u\|_Q \| v \| $$
for arbitrary $v\in \dom Q$. Since the latter is dense in $L^2$, we obtain the estimate
$$
\|  \square(\chi u) -\chi\square u\|\lesssim
\|\chi\|_{W^{2,\infty}}\|u\|_{\square} .
$$
Since
$$
\|\chi\square u\| \lesssim \|\chi\|_{W^{2,\infty}} \|u\|_{\square}
$$
is obvious, we arrive at the desired estimate. Since $\domo\square$ is a core
for $\square$ the assertion carries over to arbitrary $u\in\dom\square$.
\end{proof}

Before going on, let us note that due to the invariance under the group action
and the compact quotient our manifold has bounded geometry. We rely on \cite{Schick-96} for the
definition and a number of nice technical properties that come with bounded
geometry. The first is the existence of $r_c>0$ so that the geodesic collar
$$j:N=[0,r_c)\times bM\to M, (t,x)\mapsto \exp_x(t\nu_x)$$
is a diffeomorphism onto its image, with $\nu_x$ denoting the unit inward
normal
vector at $x$; so $t$ refers to the distance $\rho$ to the boundary mentioned
previously. Denote $j([0,\frac13 r_c)\times bM)=:N_{\frac13}$ and define
$N_{\frac23}$ accordingly.

\begin{lemma}\label{fundamental} Let $U\subset N_{\frac23}$ be a special
boundary chart and $\varphi\in
C_c^{\infty}(U,\RR)$. Then there is $\psi\in
C_c^{\infty}(U,\RR)$ so that
\begin{equation}\label{cond2} \psi|_{bM} = \varphi|_{bM}, \quad \bar L_n
\psi|_{bM} = 0.\end{equation}
Moreover, if $\varphi|_V=1$ on a set of the form $V=j([0,r)\times
R)\subset U$ with $R$ a relatively open subset of $bM$, then $\psi|_V=1$.
\end{lemma}
\begin{proof} We set $\zeta=\psi-\varphi$, so we want the derivatives
of $\zeta$ to satisfy
\[-iJ(\bar L_n - L_n)[\varphi + \zeta]|_{bM}  = (\bar L_n - L_n)[\varphi +
\zeta]|_{bM}.\]
We should have $\zeta|_{bM} = 0$, so we obtain $(\bar L_n -
L_n)\zeta|_{bM} = 0$ since the vector field is tangent to $bM$, thus
\[-iJ(\bar L_n - L_n)[\varphi + \zeta]|_{bM}   = (\bar L_n -
L_n)\varphi|_{bM}\]
\noindent
and so
\begin{align}-iJ(\bar L_n - L_n)\zeta|_{bM}   = &(\bar L_n -
L_n)\varphi|_{bM} +iJ(\bar L_n - L_n)\varphi|_{bM}\notag\\
=& 2\bar L_n\varphi|_{bM}\notag.\end{align}
Following the computations in \cite[\S3.2]{PS}, one can derive that $L_n\zeta =
d\zeta(L_n) = \langle d\zeta,\omega^n\rangle_\Lambda$ and likewise $\bar
L_n\zeta = \langle d\zeta,\bar\omega^n\rangle_\Lambda$, so that
\[iJ(\bar L_n- L_n)\zeta = (L_n+\bar L_n)\zeta = \sqrt 2\langle
d\zeta,d\rho\rangle_{\Lambda^1} = \sqrt
2\frac{\partial\zeta}{\partial\rho}.\]
\noindent
In the special boundary chart $U$, we are left with solving the equations
\[
\left\{
\begin{array}{ccc}
\left.\frac{\partial\zeta}{\partial\rho}(t,\rho)\right|_{\rho=0} &  =  &  -\sqrt 2 \bar
L_n\varphi|_{bM}   \\
 \, \zeta(t,0)\qquad & = &  0.
\end{array}
\right.
\]
\noindent
Define now, for $r<\frac23 r_c$,
\[\zeta(t,r):=-\sqrt 2\int_0^r d\rho\, \bar L_n\varphi(t,\rho).\]
 It follows that a solution $\psi$ to Eq.\ \eqref{cond2} exists. Clearly, it
satisfies the required bound on the derivatives as well as the assertion on
the level sets.
\end{proof}
\begin{prop}\label{goodcuts} There exists a good cutoff-exhaustion of $M$.
\end{prop}
\begin{proof}
We begin by constructing  a sequence of functions with bounded derivatives that
converge to $1$. To this end, let $(\varphi_i)_{i\in\ZZ}$ be a partition of
unity as in \cite[Lemma 3.22]{Schick-96}. with supports of diameter smaller
than $\frac13 r_c$. Moreover, there is a uniform bound on the number of $j$'s
so that the support of $\varphi_j$ meet a given point.
We fix $x_0\in M$ and let
$$
I^{(0)}_k:=\{ i\in\ZZ\mid \supp\varphi_i\cap B_{k+r_c}(x_0)\cap
N_{\frac13}\not=\emptyset\},
$$
Note that $\varphi^{(0)}_k:=\sum_{i\in I^{(0)}_k}\varphi_i$  satisfies
$$
1\ge \varphi^{(0)}_k\ge 1_{B_{k+r_c}(x_0)\cap N_{\frac13}} .
$$
Due to the uniform bounds for the partition of
unity,
$$
\sup_k\| \varphi^{(0)}_k\|_{W^{m,\infty}}<\infty .
$$
Note that $\varphi^{(0)}_k\in C^{\infty}_c(U,\RR)$, where $U$ is the interior
(in $\overline{M}$) of $B_{k+r_c+1}(x_0)\cap N_{\frac23}$.
The functions $\varphi^{(0)}_k$ build the ``boundary part'' of a smooth
exhaustion we want to construct. We will now modify them in a way to make sure
that the product with any function in $\domo$ in the domain $\dom(\square)$.

To this end we use Lemma \ref{fundamental} to find $\psi^{(0)}_k$ for
$\varphi^{(0)}_k$ so
that $\psi^{(0)}_k$ satisfies the requirement from (\ref{cond2}), {\it mutatis
mutandis}. Moreover,
$$
\psi^{(0)}_k|_{[0,\frac13 r_c)\times (B_k(x_0)\cap bM)} =1 ,
$$
since $\varphi^{(0)}_k$ is $1$ on the respective set by definition and the
triangle inequality.

Denote $I_k:=\{ i\in\ZZ\mid \supp\varphi_i\subset
B_{k+r_c}(x_0)\}$ and
$$
\psi^{(1)}_k:= (1-\psi^{(0)}_k)\cdot \sum_{i\in I_k}\varphi_i .
$$
By the assumption on the support of the $\varphi_i$, the sum is $1$ on
$B_k(x_0)$ and $\supp\psi^{(1)}_k\subset
B_{k+r_c}(x_0)$. In particular, $\psi^{(1)}_k\in C_c^\infty(M,\RR)$ and
$$
\chi_k:= \psi^{(0)}_k+\psi^{(1)}_k
$$
is $1$ on $B_k(x_0)$. Thus (C1) from the definition of a good cutoff
exhaustion above is satisfied. The function $\psi^{(0)}_k$ was constructed so
that the required condition (\ref{tazera}) holds and since $\psi^{(1)}_k$ is
supported away from the boundary, $\chi_k$ satisfies (C2). The uniform bound on
the derivatives is evident from the definition and the properties of the
partition of unity.
\end{proof}

\begin{proof}[Proof of Theorem 3]
 We have to show that any $u\in\dom\square$ can be approximated by a sequence
$(u_k)$ in $\ccpq$ in the Hilbert space
$(\dom\square,\|\cdot\|_\square)$. Since $\domo\square$ is dense by Prop.\ \ref{dbarlap} we can restrict to the case in which $u\in \domo\square$. Since
$\ccpq$ is convex, its weak and norm closures in
$(\dom\square,\|\cdot\|_\square)$ coincide, so we are left with finding $(u_k)$
that converges weakly in the latter space. We take a good cutoff exhaustion
$(\chi_k)$ and claim that $u_k:=\chi_k u$ does the job.
By (C1) and (C3) we know that $u_k\to u$ in $\lzpq$ as $k\to\infty$. Moreover,
by (\ref{chi_u}) in Prop.\ \ref{prop7.4} it follows that $(u_k)$ is
bounded in $(\dom\square,\|\cdot\|_\square)$. It thus has a weakly convergent
subsequence that has to converge to $u$ by uniqueness of the limit and the
$L^2$-convergence we  already established.\end{proof}

%%%%%%%%%%%%%%%%%%%%%%%%%%%%%%%%
\begin{ack} It is a pleasure to thank Emil Straube for his help with the cutoff
argument in Sect.\ \ref{essential}. JJP would like to thank the Math Department
at TU-Chemnitz and the Erwin Schr\"odinger
Institute for extremely fruitful visits.\end{ack}
%%%%%%%%%%%%%%%%%%%%%%%%%%%%%%%%
%%%%%%%%%%%%%%%%%%%%%%%%%%%%%%%%
%%%%%%%%%%%%%%%%%%%%%%%%%%%%%%%%

\end{document}